\documentclass[11pt]{article}
\usepackage{amscd}
\usepackage{verbatim}
\usepackage{amssymb}
\usepackage{amsmath}
\usepackage{amsthm}
\usepackage{graphics,epsfig}

\begin{document}
\vskip 0 true cm \flushbottom

\begin{center}
\vspace{24pt} { \large \bf The Ricci flow of the ${\mathbb
R}{\mathbb P}^3$-geon and
noncompact manifolds with essential minimal spheres} \\
\vspace{30pt}
{\bf T Balehowsky}\footnote{balehows@ualberta.ca} and
{\bf E Woolgar}\footnote{ewoolgar@math.ualberta.ca} 

\vspace{24pt} 
{\footnotesize $^\dag$ Dept of Mathematical and Statistical
Sciences,
University of Alberta,\\
Edmonton, AB, Canada T6G 2G1.}
\end{center}
\date{\today}
\bigskip

\begin{center}
{\bf Abstract}
\end{center}
\noindent It is well-known that the Ricci flow of a closed
3-manifold containing an essential minimal 2-sphere will fail to
exist after a finite time. Conversely, the Ricci flow of a complete,
rotationally symmetric, asymptotically flat manifold containing no
minimal spheres is immortal. We discuss an intermediate case, that
of a complete, noncompact manifold with essential minimal
hypersphere. For 3-manifolds, if the scalar curvature vanishes on
asymptotic ends and is bounded below initially by a negative
constant that depends on the area of the minimal sphere, we show
that a singularity develops in finite time. In particular, this
result applies to asymptotically flat manifolds, which are a
boundary case with respect to the neckpinch theorem of M Simon. We
provide numerical evolutions to explore the case where the initial
scalar curvature is less than the bound.


\setcounter{equation}{0}
\newpage

\section{Introduction}
\setcounter{equation}{0}

\noindent A Ricci flow on a manifold $M$ is a family of
Riemannian metrics $g_{ij}(t;x)$, $x\in M$, $t\in I\subseteq
{\mathbb R}$ the family parameter, and $I$ a connected interval,
satisfying
\begin{equation}
\frac{\partial g_{ij}}{\partial t} = -2R_{ij}\ , \label{eq1.1}
\end{equation}
where $R_{ij}$ is the Ricci tensor of $g_{ij}(t,x)$. It is usually
more convenient to study the {\it Hamilton-DeTurck flow} (or {\it
Ricci-DeTurck flow})
\begin{equation}
\frac{\partial g_{ij}}{\partial t} = -2R_{ij}+\pounds_X g_{ij}
\ . \label{eq1.2}
\end{equation}
If $g$ solves (\ref{eq1.2}), then the pullback $\psi_t^*g$ solves
(\ref{eq1.1}) (\cite{CK}, p 80), where $\psi_t$ is a family of
time-dependent diffeomorphisms generated by the vector field $X$.
From (\ref{eq1.2}), the scalar curvature $R$ of $g_{ij}(t)$ evolves
according to
\begin{equation}
\frac{\partial R}{\partial t} = \Delta R + \nabla_X R +2R_{ij}R^{ij}
\ , \label{eq1.3}
\end{equation}
where $\Delta:=g^{ij}\nabla_i \nabla_j$ is the ($t$-dependent)
Laplacian.

Ricci flow of asymptotically flat manifolds arises in several
physical and mathematical contexts, ranging from the physics of
closed string tachyon condensation \cite{GHMS} to existence problems
for static Einstein metrics \cite{GOW}. It is known that
asymptotically flat initial data remain asymptotically flat and
smooth when evolved by Ricci flow for some time interval $[0,T)$
\cite{DM, OW}. For complete, rotationally symmetric initial data
with no minimal sphere present (corresponding in static general
relativity to the absence of black hole horizons), it has been shown
that the flow exists for all future time and converges to flat space
\cite{OW}.

Then the question arises as to what happens when complete,
rotationally symmetric initial data containing a minimal sphere are
evolved. There are two cases, depending on whether the minimal
sphere is topologically essential (representing a nontrivial class
in $\pi_2(M)$) or, as will be the case for data on ${\mathbb R}^n$,
inessential. Of these, the essential case is the easier one to
study, and is the subject of the present paper.

It is well-known that a closed 3-manifold admitting an essential
minimal 2-sphere, when evolved under the Ricci flow, will develop a
singularity within finite time. Depending on the initial
configuration, the singularity may be localized at or near the
minimal sphere, or may be global in the sense that the manifold
collapses everywhere at that time. It is expected that this result
carries over to the noncompact case. Somewhat to the contrary
though, asymptotically flat initial data are critical data with
respect to the pinching theorem of M Simon \cite{Simon}.  That is,
Simon proves that warped products of a line with a positively curved
closed manifold will evolve to form a neckpinch singularity in
finite time, provided the initial data obey certain asymptotic
conditions at the ends of the line, including a condition on
asymptotic growth of the area of the closed manifold factor. Simon's
growth condition is written as a strict (i.e., open) inequality
which is not satisfied by asymptotically flat manifolds, but such
manifolds would lie in the closure. This
suggests that such data may exhibit interesting evolutions,
including perhaps critical phenomena of the sort observed in certain
numerical evolutions in general relativity in a scenario sometimes
called ``critical collapse'' \cite{Choptuik}.

Garfinkle and Isenberg \cite{GI1} studied numerical Hamilton-DeTurck
flow of a 3-sphere with ``corsetted'' initial metric, admitting an
inessential minimal 2-sphere (the ``waist''). For tight corsetting,
meaning that the waist has very small area relative to the $2/3^{\rm
rds}$-power of the volume of the 3-sphere, the waist ``pinches off''
(a local singularity forms there), whereas for more gentle
corsetting, the entire sphere shrinks to a point before a local singularity can form. The critical solution separating these two alternatives is a
degenerate neckpinch singularity modelled by the Bryant soliton
\cite{GI2}.

The Garfinkle-Isenberg result shows that initial data for the Ricci
flow with an inessential minimal sphere divides into two disjoint
sets whose common boundary consists of points whose evolutions
exhibit critical behaviour.

Husain and Seahra \cite{HS} then considered a numerical
Hamilton-DeTurck flow of a sequence of initial metrics, each
rotationally symmetric and reflection symmetric through an essential
minimal 2-sphere, sometimes referred to as the ``throat'' or
``bridge''. Evolution occurred on a bounded region, with boundary
conditions imposed. Two alternatives were again found; for some
initial data, the throat pinched off while, for other data, the
throat expanded to infinity.

We study a related question. Consider simply connected, {\em
noncompact, complete} manifolds with rotational symmetry. These are
${\mathbb R}^n$, ${\mathbb R}\times S^{n-1}$, and certain quotients
thereof. As ${\mathbb R}^n$ has no essential minimal sphere,
consider ${\mathbb R}\times S^{n-1}$, and endow it with an ${\rm
SO}(n)$-symmetric metric such that there is a minimal hypersphere
located at, say, $r=0$ and an isometry corresponding to reflection
in that hypersphere. One such metric is the $t=0$ slice of the
$(n+1)$-dimensional Schwarzschild-Tangherlini metric. One can then
identify points under the action ${\mathbb R}\times
S^{n-1}\ni(r,p)\mapsto (-r,-p)$, where $p\in S^{n-1}$ and $-p$ is
the antipode of $p$. This produces a smooth metric on ${\mathbb
R}\times{\mathbb R}{\mathbb P}^{n-1}\cong{\mathbb R}{\mathbb P}^{n}
\setminus \{ {\rm pt} \}$. For $n=3$, this model is known in
gravitational physics as the {\it ${\mathbb R}{\mathbb P}^{3}$
geon}. We pose and will answer the question, ``What is the Ricci
flow evolution of the ${\mathbb R}{\mathbb P}^{3}$ geon?''\footnote
{This was posed by DM Witt to the second author quite some time
ago.}

In section 2, we adapt a standard Ricci flow argument from the
setting of closed manifolds to that of asymptotically flat
manifolds. This forms the basis for what follows, and shows that
essential 2-spheres in asymptotically flat 3-manifolds collapse
whenever the initial scalar curvature is bounded below by a
nonpositive constant that can depend on the initial area of the
minimal sphere. The result is valid for arbitrary Hamilton-DeTurck
flow, including Ricci flow. As the covering space of the ${\mathbb
R}{\mathbb P}^{3}$ geon obeys this bound, this suffices to answer
the question just posed.

This raises the possibility that, by choosing initial data for which
the lower bound on initial scalar curvature is violated (e.g., by
choosing a different initial metric on the geon manifold, with the
same isometries), interesting dynamics might arise, such as observed
in the studies cited above \cite{GI1, GI2, HS}. It also raises the
question of the precise comparison, if any, of this result to the
numerical work of \cite{HS}, who found that for certain initial data
with scalar curvature well above our lower bound, their numerical
Hamilton-DeTurck evolution did not lead to collapse.

To understand these issues, we perform our own numerical
simulations. In section 3, we lay the groundwork. We first discuss
Hamilton-DeTurck flow with the ``DeTurck trick'' formulation which
we use for our numerical evolutions in Section 4 (and which
\cite{GI1} used). For comparison purposes, we then discuss the
normal coordinate Hamilton-DeTurck formulation used in \cite{HS}.

Section 4 contains our numerical results. We use the same form of
initial data as Husain and Seahra \cite{HS} and, for practical
purposes, we also now restrict our evolution to a bounded manifold,
but we use different evolution equations\footnote
{We use a different Hamilton-DeTurck system. Occurrence of collapse of the minimal surface should not depend on this choice; see section 2.}
and different boundary conditions at the
boundary at large $r$. Our numerical evolutions always exhibit
collapse of the throat. Section 5 contains a brief discussion of the
numerical evolution.

Our convention for the curvature tensor $R^a{}_{bcd}$ is that used
in \cite{HE} and equals the quantity denoted by $R_{cdb}{}^a$ in
\cite{CK}. We write $R_{abcd}:=g_{ae}R^e{}_{bcd}$. We denote the
Laplacian by $\Delta:=g^{ab}\nabla_a\nabla_b$. For a definition of
asymptotic flatness, see \cite{OW}.

\section{Time derivative of the area of a minimal sphere}
\setcounter{equation}{0}

\subsection{The maximum principle for scalar curvature}

\noindent For Ricci flow on compact manifolds, it is standard that
positive scalar curvature is associated to contraction and
concentration of curvature under the Ricci flow. In this subsection
is that scalar curvature, we recall that this is also true for
complete manifolds with vanishing scalar curvature at infinity.

\bigskip

\noindent {\bf Proposition 2.1.} {\sl Assume that (\ref{eq1.3}) has
a solution on some time interval $[0,T]$, $T>0$, on a complete
manifold $M$ with one or more asymptotic ends, and assume that
$R(t,x)\to 0$ as $x\to\infty$ (that is, as the point $x$ tends to an
asymptotic end), for all $t\in [0,T]$. Let the scalar curvature
$R_0(x)$ of the initial metric $g_0(x):=g(0,x)$ obey
\begin{equation}
\inf_{x\in {\cal M}} \left \{ R_0(x) \right \}=:-a^2\ , \label{eq2.1}
\end{equation}
for some $a\ge 0$. Then}
\begin{equation}
R(t,x)\ge \frac{-a^2}{1+\frac{2a^2}{n}t}\ . \label{eq2.2}
\end{equation}

\bigskip

\noindent{\bf Remark 2.2.} The existence assumption always holds for
Ricci flow developing from asymptotically flat initial data
(\cite{OW}, \cite{DM}).

\bigskip

\noindent{\bf Proof.} Equation (\ref{eq1.3}) can be written as
\begin{equation}
\frac{\partial R}{\partial t} = \Delta R + \frac{2}{n}R^2
+ 2\left \vert R_{ij}-\frac{1}{n}g_{ij}R \right \vert^2
+\nabla_X R\ . \label{eq2.3}
\end{equation}
Exhaust $M$ by a sequence of closed, bounded sets $K_i$. An easy
application of the maximum principle (e.g., \cite{Jost}) to the
$i^{\rm th}$ subset yields that the minimum of $R$ occurs on the
parabolic boundary of $[0,T]\times K_i$. Now take $i\to\infty$ and
use that $R(t,x)\to 0$ as $x\to\infty$ to deduce that the minimum
must either be zero or must occur on the initial data. This proves
that
\begin{equation}
R(t,x)\ge \inf_{x\in {\cal M}} \left \{ R_0(x) \right \}=-a^2,\ a\ge 0.
\label{eq2.4}
\end{equation}
If $a=0$, this proves the theorem, so we now consider $a>0$, in which case the infimum is a minimum.

For $a>0$, let $Q := -\frac{1}{a^2} \left ( 1 +
\frac{2a^2}{n}t \right ) R$. Since $R=-a^2<0$ at its minimum, then
$Q>0$ there. Hence the maximum of $Q$ is positive, and it then
follows from the definition of $Q$ that, at the point where $Q$
achieves its maximum, $R$ is negative. As well, $Q(t,x)\to 0$ as
$x\to \infty$.

Work on the intervalcompact set $[0,T]\times K_i$. From
(\ref{eq2.3}), we have
\begin{equation}
\frac{\partial Q}{\partial t} \le \Delta Q +\nabla_X Q
+ \frac{2}{n} R \left ( Q-1 \right )\ . \label{eq2.5}
\end{equation}
Then from the maximum principle (and since $R<0$ where $Q$ achieves
its maximum), either the maximum of $Q$ occurs on the parabolic
boundary of $[0,T]\times K_i$ or $Q-1\le 0$ at the maximum. If it
occurs on the parabolic boundary, by taking $i\to\infty$ and using
that $Q(t,x)\to 0$ as $x\to\infty$, we see that the maximum must
occur on the initial boundary, where $Q(0,x)=-\frac{1}{a^2}R(0,x)\le
1$. Hence, $Q(t,x)\le 1$, and the result follows. \qed

\bigskip

\noindent{\bf Remark 2.3.} Obviously, these results also hold on
$[0,T]\times K$, $K\subset M$ a compact set with boundary $\partial
K$, provided that $R(t,x)=0$ for $x\in\partial K$. But they do not
hold if merely $R(t,x)\ge -const$ on $\partial K$, an observation
which appears to be relevant in numerical studies (see section
3.2.2).

\bigskip

\subsection{Evolution of hypersurfaces}

\noindent In this subsection, we obtain the evolution equation for
the area of a closed hypersurface $\Sigma$ in an $n$-manifold
$({\cal M}^n,g)$, defined (locally at least) by an expression of the
form $F(x^i)=0$, where $x^i$ are the coordinates held constant in
the derivative $\frac{\partial}{\partial t}$ appearing in
(\ref{eq1.2}).\footnote
{For hypersurfaces defined more generally by an equation of the form
$F(t,x^i)=0$, the generalization of our result is easy to obtain,
using the transport theorem.}

Let $n^a$ be a smooth unit vector field normal to $\Sigma$. Let
$h_{ab}$ be the induced metric and let $H_{ab}$ be the extrinsic
curvature of $\Sigma$, with trace $H$ (thus $H$ is the mean
curvature, taken here to mean the sum rather than the average of the
principal curvatures at a point). If ${\tilde H}_{ab}$ denotes the
trace-free part of $H_{ab}$, then
\begin{equation}
H_{ab}={\tilde H}_{ab}+\frac{1}{n-1}h_{ab}H\ . \label{eq2.6}
\end{equation}
If ${\cal R}$ denotes the scalar curvature of the induced metric
$h_{ab}$ on $\Sigma$, then
\begin{equation}
R - {\rm Ric}(n,n)= \frac12 \left [ {\cal R}
- \left ( \frac{n-2}{n-1}\right )
H^2 +R+\vert {\tilde H} \vert^2 \right ] \label{eq2.7}
\end{equation}
by the Gauss equation. Finally, let the induced area element on
$\Sigma$ be $d\mu$. From the Hamilton-DeTurck flow equation
(\ref{eq1.2}), the area $\vert \Sigma \vert:=\int_{\Sigma}d\mu$ of
$\Sigma$ evolves as
\begin{eqnarray}
\frac{d}{dt}\vert \Sigma \vert&=&\int_{\Sigma}\frac{\partial}{\partial t}d\mu
=\frac12\int_{\Sigma}h^{ij}\frac{\partial h_{ij}}{\partial t}d\mu
=\frac12\int_{\Sigma}h^{ij}\frac{\partial g_{ij}}{\partial t}d\mu
\nonumber\\
&=&-\int_{\Sigma}h^{ij}\left ( R_{ij}-\nabla_i X_j \right )d\mu
=-\int_{\Sigma}\left ( R-R_{ij}n^in^j-h^{ij}\nabla_i X_j\right ) d\mu
\nonumber\\
&=&-\frac12 \int_{\Sigma} \left ( {\cal R} - \left ( \frac{n-2}{n-1}\right )
H^2 +R+\vert {\tilde H} \vert^2-2Hn\cdot X\right )  d\mu
\ ,\label{eq2.8}
\end{eqnarray}
using (\ref{eq2.7}), the divergence theorem on $\Sigma$, and
$H:=h^{ij}\nabla_i n_j$.

\subsection{3-manifolds}

\noindent Now fix the dimension to be $n=3$. Then we have
$\int_{\Sigma} {\cal R} d\mu=4\pi\chi(\Sigma)$, where $\chi(\Sigma)$
is the Euler characteristic of $\Sigma$. Then
\begin{equation}
\frac{d}{dt}\vert \Sigma \vert=-2\pi\chi(\Sigma)+ W(\Sigma)
-\frac12 \int_{\Sigma} \left ( R+\vert {\tilde H} \vert^2 \right ) d\mu
+ \int_{\Sigma} Hn\cdot X d\mu\ , \label{eq2.9}
\end{equation}
where $W(\Sigma):=\frac14 \int_{\Sigma} H^2 d\mu$ is the {\it
Willmore energy} \cite{Willmore} of $\Sigma$.\footnote
{When $\Sigma$ is a 2-sphere, the combination of Willmore energy and
$\chi(\Sigma)$ that appears in (\ref{eq2.9}) also appears in the
definition of the Hawking quasi-local mass $m_H:=\frac{\vert \Sigma
\vert^{1/2}}{16\pi^{3/2}} \left [ 4\pi-W(\Sigma)\right ]$.}
Using Proposition\ 2.1, (\ref{eq2.9}) becomes
\begin{equation}
\frac{d}{dt}\vert \Sigma \vert \le -2\pi\chi(\Sigma) + W(\Sigma)
+ \frac{3a^2\vert \Sigma \vert}{6+4t}
+ \int_{\Sigma}Hn\cdot X d\mu\ . \label{eq2.10}
\end{equation}
The tightest bound occurs in the case of a minimal 2-sphere. Then
$H\vert_{\Sigma}$ and the Willmore energy vanish and $\chi(S^2)=2$
so
\begin{equation}
\frac{d}{dt}\vert \Sigma \vert \le -4\pi
+ \frac{3a^2}{6+4t}\vert \Sigma \vert \ . \label{eq2.11}
\end{equation}
We note that the initial time derivative of $|\Sigma|$ will be
negative whenever
\begin{equation}
a^2<\frac{8\pi}{|\Sigma|_0}\ , \label{eq2.12}
\end{equation}
where $|\Sigma|_0$ denotes the initial area of $\Sigma$. Even when
not initially negative, the derivative can turn negative and remain
so thereafter.

Now we can regard the left-hand side of (\ref{eq2.11}) as the
forward difference quotient of the functional whose value is the
area of the smallest essential minimal surface present at $t$. Then
(see Lemma 2.22 of \cite{MT}) an upper barrier for the area of the
smallest minimal surface at any $t\ge 0$ is provided by the function
$\Psi(t)$ that solves the initial value problem
\begin{eqnarray}
\frac{d\Psi}{dt}-\frac{3a^2\Psi}{6+4t} &=& -4\pi\ , \label{eq2.13}\\
\Psi(0)=\vert \Sigma\vert_0\ . \nonumber
\end{eqnarray}
The solution is
\begin{equation}
\Psi(t) = \begin{cases} \frac{24\pi}{3a^2-4}(1+\frac23 t)
+\left ( |\Sigma|_0 - \frac{24\pi}{3a^2-4} \right )
\left ( 1+\frac23 t\right )^{3a^2/4}  & \text{for $a^2\neq 4/3$,}\\
(1+\frac23 t)\left [ |\Sigma|_0 - 6\pi \log\left ( 1+\frac23 t\right )
\right ] & \text{for $a^2= 4/3$.}  \label{eq2.14}
\end{cases}
\end{equation}

\bigskip
\noindent {\bf Proposition 2.4.} {\sl Say that the conditions of
Proposition 2.1 hold and that $(M,g_0)$ contains a minimal sphere of
area $4\pi\delta$ for some $\delta>0$. Let $3a^2-4<6/\delta$. Then
the flow fails to exist in finite time.}

\bigskip
\noindent{\bf Proof.} We must show that if the flow lasts long
enough, the minimal sphere collapses. If $a^2=4/3$, this is
immediate from the second line of (\ref{eq2.14}), since $\Psi$ has a
zero at some $t>0$. Otherwise, let $\gamma:= 3a^2-4\neq 0$. The top
line of (\ref{eq2.14}) becomes
\begin{equation}
\Psi(t)= 4\pi \left ( 1+\frac23 t \right ) \left [ \frac{6}{\gamma}
+ \left (\delta - \frac{6}{\gamma} \right ) \left ( 1 +\frac23 t
\right )^{\gamma/4} \right ] \ , \label{eq2.15}
\end{equation}
so $\Psi(t)=0$ for some $t>0$ iff
\begin{equation}
\left ( 1+\frac23 t \right )^{\gamma/4}=\frac{6}{6-\gamma\delta}\ . \label{eq2.16}
\end{equation}
Now $\frac{6}{6-\gamma\delta}>1$ if $0<\gamma<6/\delta$ and
$0<\frac{6}{6-\gamma\delta}<1$ if $\gamma<0$. The left-hand side
equals 1 at $t=0$ and is otherwise monotonic, increasing without bound if
$\gamma>0$, and decreasing with asymptote $0$ if $\gamma<0$. Hence there is always a root for $t>0$ if $\gamma<6/\delta$.\qed

\bigskip

\noindent{\bf Corollary 2.5.} For any $X$, consider the flow
(\ref{eq1.3}) developing from an initial metric which describes the
${\mathbb R}{\mathbb P}^3$-geon. Then the flow fails to exist in
finite time.

\bigskip
\noindent{\bf Proof.} The Riemannian double cover of the geon is the
$\delta=1$, $a=0$ case above. Furthermore, since the initial data
are asymptotically flat, then $R(t,r)\to 0$ as $r\to\infty$ for all
$t$ during the flow (\cite{DM, OW}). Then the result holds on the
Riemannian double cover and, hence, on the manifold itself because
the covering map is a local isometry.

Alternatively, the geon itself has $\chi(M)=1$ and $|\Sigma|=2\pi$,
so equation (2.13) applies, with the $4\pi$ on the right-hand side
of the differential equation replaced by $2\pi$ and
$|\Sigma_0|=2\pi$. The solution (\ref{eq2.14}) is then multiplied by
an overall factor of $1/2$, which has no effect in the
proposition.\qed

\bigskip

\subsection{2-dimensions}

\noindent The $n=2$ case is itself of some interest. Then
(\ref{eq2.8}) becomes
\begin{equation}
\frac{d}{dt}\vert \Sigma \vert = -\frac12 \int_{\Sigma} \left ( R
-2Hn\cdot X  \right )  d\mu\ , \label{eq2.17}
\end{equation}
with $H$ the geodesic curvature of the closed curve $\Sigma$. Then
\begin{equation}
\frac{d}{dt}\vert \Sigma \vert = -\frac12 \int_{\Sigma} Rd\mu \label{eq2.18}
\end{equation}
if $\Sigma$ is a geodesic. If $R\to 0$ at infinity (or on any
boundary), Proposition 2.1 yields that $R\ge -\frac{a^2}{1+a^2t}$,
and so (2.19) becomes
\begin{equation}
\frac{d}{dt}\vert \Sigma \vert \le \frac{a^2}{2(1+a^2t)}\vert \Sigma \vert
\ . \label{eq2.19}
\end{equation}
Thus
\begin{equation}
\vert \Sigma \vert(t) \le \vert \Sigma \vert (0) \sqrt{1+a^2t}\ . \label{eq2.20}
\end{equation}

\section{Rotational symmetry and Hamilton-DeTurck flow}
\setcounter{equation}{0}

\noindent There are two scales in this problem, the initial area of
the minimal surface and the minimum of the initial scalar curvature.\footnote
{This presumes the scalar curvature has a negative minimum. In the initial data we study, the only exception will be the Schwarzschild-Tangherlini data, for which $R(0,p)=0$ for all points $p$.}
The metric can
be rescaled to fix one but not both of these. By choosing
$\delta=1$, we fix the initial area. There is then the intriguing
possibility that collapse of the minimal sphere will not occur if
the geon initial data is replaced with initial data with scalar
curvature sufficiently negative that $a^2>10/3$. In the next section
we will study this question numerically. In this section, we
formulate the equations we will need.

There are several versions of Hamilton-DeTurck system that are used
in the study of rotationally symmetric Ricci flow. We discuss the
formulation we will use for our numerical integrations in the first
subsection. We then discuss, for comparison purposes, other
formulations that have been used in related work.

\subsection{DeTurck's background connection method}

\subsubsection{The flow equations with a background connection}

\noindent The first Hamilton-DeTurck system arises from DeTurck's
original trick for proving short-time existence of Ricci flow, and
results in a parabolic strictly system of two equations. The idea is
that parabolicity fails only because there are families of flows
whose members are distinguished from each other only by a continuous
time-dependent deformation of the coordinates. Thus, parabolicity is
restored by ``breaking coordinate invariance''. This is done by
fixing a $t$-independent background connection throughout the flow.
That is, let $\Gamma^i_{jk}$ be the flowing connection (compatible
with $g_{ij}(t)$) written in local coordinates, and let
$\breve{\Gamma}^i_{jk}$ denote the chosen background connection in
these coordinates. Define the vector field
\begin{equation}
X^i:=g^{jk}\left ( \Gamma^i_{ij}-\breve{\Gamma}^i_{jk} \right )
\ . \label{eq3.1}
\end{equation}
Then the system (\ref{eq1.2}) is parabolic (for any choice of $t$-independent
connection $\breve{\Gamma}$). For rotationally symmetric flow on
${\mathbb R}^n$, an obvious choice of coordinates is
\begin{equation}
ds^2 = e^{2A(t,r)} dr^2 + r^2 e^{2B(t,r)}g(S^{n-1},{\rm can})\ ,
\label{eq3.2}
\end{equation}
where $g(S^{n-1},{\rm can})$ denotes the metric of unit sectional
curvature on $S^{n-1}$. Then it makes sense to choose a background
connection to arise from a rotationally symmetric background metric.
Since every rotationally symmetric metric is conformal to a
flat metric, we write the background metric as
\begin{equation}
\breve{g}=e^{2\psi(r)}\left [ dr^2 +r^2 g(S^{n-1},{\rm can})
\right ]\ , \label{eq3.3}
\end{equation}
for some function $\psi(r)$ that must be specified. Choices of
$\psi$ with rotational symmetry and asymptotic flatness include
\begin{equation}
\psi = \begin{cases} 0 \ , & \text{flat background,} \\
  \frac{2}{\alpha}\ln \left [ \left ( 1+
  r^{-\alpha}\right ) \right ] \ , & \alpha=const>0. \label{eq3.4}
\end{cases}
\end{equation}
The last choice listed is a family of metrics which were used as
initial data for Ricci flow (with $1\le\alpha\le 2$ and $n=3$ dimensions) in \cite{HS}
and which would, if $\alpha=n-2$, correspond to a time-symmetric
slice of $(n+1)$-dimensional Schwarzschild spacetime.

Then the DeTurck vector field is
\begin{eqnarray}
X=:V\frac{\partial}{\partial r} &=&\bigg \{ e^{-2A} \left [
\frac{\partial A}{\partial r} - \psi'(r) \right ]\label{eq3.5}\\
&&-(n-1) \left [ \frac{e^{-2A}-e^{-2B}}{r}
+e^{-2A} \frac{\partial B}{\partial r}
-e^{-2B}\psi'(r)\right ] \bigg \}\frac{\partial}{\partial r}\ . \nonumber
\end{eqnarray}
The Hamilton-DeTurck system becomes
\begin{eqnarray}
\frac{\partial A}{\partial t} &=& e^{-2A}\left \{
\frac{\partial^2 A}{\partial r^2} - \left (
\frac{\partial A}{\partial r} \right )^2 +(n-1)
\left ( \frac{\partial B}{\partial r} +\frac{1}{r} \right )^2-\psi''(r)
+\frac{\partial A}{\partial r}\psi'(r) \right \}
\nonumber\\
&&+(n-1)e^{-2B} \left \{ \frac{1}{r}\frac{\partial A}{\partial r}
-\frac{2}{r}\frac{\partial B}{\partial r}-\frac{1}{r^2} +\psi''(r)
+\frac{\partial A}{\partial r}\psi'(r)
-2\frac{\partial B}{\partial r}\psi'(r) \right \}\ ,\label{eq3.6}\\
\frac{\partial B}{\partial t} &=& e^{-2A}\left \{
\frac{\partial^2 B}{\partial r^2} -\frac{1}{r^2}
-\left ( \frac{\partial B}{\partial r}+\frac{1}{r}\right )
\psi'(r)\right \}\nonumber \\
&&+(n-1)e^{-2B}\left ( \frac{\partial B}{\partial r}+\frac{1}{r}\right )
\left ( \psi'(r)+\frac{1}{r}\right )-\frac{(n-2)}{r^2}e^{-2B}
\ .\label{eq3.7}
\end{eqnarray}

If we define
\begin{equation}
S:=V_{n-1}r^{n-1}e^{(n-1)B}\label{eq3.8}
\end{equation}
to be the area of $r=const$ spheres, with $V_{n-1}$ being the
volume of an $(n-1)$-sphere
with unit sectional curvature, we can rewrite the above system
with $B$ replaced by $S$.

In $n=3$ dimensions, this gives
\begin{eqnarray}
\frac{\partial A}{\partial t} &=& e^{-2A}\left \{
\frac{\partial^2 A}{\partial r^2} - \left (
\frac{\partial A}{\partial r} \right )^2 +\frac{1}{2S^2}
\left ( \frac{\partial S}{\partial r}\right )^2-\psi''(r)
+\frac{\partial A}{\partial r}\psi'(r) \right \}
\nonumber\\
&&+\frac{8\pi r^2}{S} \bigg \{ \frac{1}{r}\frac{\partial A}{\partial r}
-\frac{1}{rS}\frac{\partial S}{\partial r}+\frac{1}{r^2}
+\psi''(r)+\left ( \frac{\partial A}{\partial r} +\frac{2}{r} \right )
\psi'(r) \nonumber\\
&&\qquad -\frac{1}{S}\frac{\partial S}{\partial r}\psi'(r) \bigg \}
\ ,\label{eq3.9}\\
\frac{\partial S}{\partial t} &=& e^{-2A}\left \{
\frac{\partial^2 S}{\partial r^2} -\frac{1}{S}\left (
\frac{\partial S}{\partial r} \right )^2
-\psi'(r)\frac{\partial S}{\partial r}\right \}
+\frac{8\pi r^2}{S}\frac{\partial S}{\partial r}
\left ( \psi'(r)+\frac{1}{r}\right ) -8\pi
\ .\label{eq3.10}
\end{eqnarray}
Once a function $\psi$ is chosen, it is clear that this system is
parabolic. This system was used in \cite{GI1} in their study of
rotationally symmetric flow on ``corsetted'' 3-spheres.

\subsubsection{Minimal surface boundary conditions with background connection}

\noindent Let us now consider the problem of the previous section,
which is the evolution of an asymptotically flat metric on a
manifold with inner boundary $\Sigma$ which is a minimal surface. At
$\Sigma$, both the mean curvature and the vector field $X$ must
vanish, so
\begin{eqnarray}
V_{\Sigma}(t)&:=&V(t,r_0)=0\ , \label{eq3.11}\\
H_{\Sigma}(t)&:=&H(t,r_0)=0\ . \label{eq3.12}
\end{eqnarray}
The second of these conditions is of course $\frac{\partial
S}{\partial r}\bigg \vert_{r_0}=0$ or, equivalently,
\begin{equation}
\frac{\partial B}{\partial r}\bigg \vert_{r_0}+\frac{1}{r_0}=0
\ , \label{eq3.13}
\end{equation}
while the first then yields
\begin{equation}
e^{-2A(t,r_0)} \left ( \frac{\partial A}{\partial r}\bigg \vert_{r_0}
-\psi'(r_0)\right ) +(n-1)e^{-2B(t,r_0)} \left (\frac{1}{r_0}+\psi'(r_0)
\right ) =0\ . \label{eq3.14}
\end{equation}

\subsection{Polar coordinate gauge}

\subsubsection{The flow equations in Gaussian polar coordinates}

\noindent This system was used by Husain and Seahra \cite{HS} and is
not parabolic. The metric is
\begin{equation}
ds^2 = dr^2+F^2(t,r)g(S^{n-1},{\rm can})\ , \label{eq3.15}
\end{equation}
corresponding to $A=0$, $B=\log (F/r)$ in (\ref{eq3.2}). The
coordinate $r$ is distance from the minimal sphere at $r=0$. In this
system, the Hamilton-DeTurck flow becomes a constrained system
comprised of the differential equation
\begin{equation}
\frac{\partial F}{\partial t} = \frac{\partial^2 F}{\partial r^2}
+ \frac{1}{F} \left ( \frac{\partial F}{\partial r}\right )^2
+V\frac{\partial F}{\partial r}-\frac{n-2}{F}\label{eq3.16}
\end{equation}
and a differential constraint
\begin{equation}
\frac{\partial V}{\partial r} = -\frac{(n-1)}{F}
\frac{\partial^2 F}{\partial r^2} \label{eq3.17}
\end{equation}
which determines the generator $X=V\frac{\partial}{\partial r}$ of
the diffeomorphism in (\ref{eq1.2}). The mean curvature of a
constant-$r$ sphere is
\begin{equation}
H:=\frac{(n-1)}{F} \frac{\partial F}{\partial r} =\frac{1}{S}
\frac{\partial S}{\partial r}\ , \label{eq3.18}
\end{equation}
where $S=V_{n-1}F^{n-1}$ is the area of an $r=const$ sphere
(cf.\ (\ref{eq3.8})), and so the solution of the constraint can be
written as
\begin{equation}
V(t,r)=C(t)+H(t,r_0)-H(t,r)-\frac{1}{(n-1)}\int\limits_{r_0}^r
H^2(t,r') dr' \ , \label{eq3.19}
\end{equation}
where $C(t)$ is an arbitrary function of $t$, to be determined by
boundary conditions at $r_0$. When the manifold contains an origin
for the rotational symmetry (e.g., ${\mathbb R}^n$), this formulation
suffers from the problem that $H$ diverges there.

In the case of $n=3$ dimensions, (\ref{eq3.16}) leads to the
remarkably simple equation
\begin{equation}
\frac{\partial S}{\partial t}
=\frac{\partial^2 S}{\partial r^2}
+V\frac{\partial S}{\partial r}-8\pi\label{eq3.20}
\end{equation}
for the area $S(t,r)$ of an orbit of the rotational symmetry.
Indeed, this equation would be linear if not for the dependence of
$V$ on $S$ via the constraint (\ref{eq3.17}) (which takes a slightly
more complicated form when $F$ is replaced by $S$).

\subsubsection{Minimal surface boundary conditions in polar coordinates}

\noindent The condition that $V$ should vanish at the boundary
$\Sigma$ is simply $C(t)=0$. The condition (\ref{eq3.12}) that the
boundary $\Sigma$ be a minimal surface is $H(t,r_0)=0$. These
conditions, applied to (\ref{eq3.19}), then imply that
\begin{equation}
V(t,r)=-H(t,r)-\frac{1}{(n-1)}\int\limits_{r_0}^r H^2(t,r')dr' .
\label{eq3.21}
\end{equation}
Since $H$ is not identically zero (for all $r$), then $V(t,r)<0$ on
every convex surface $r>r_0$. Husain-Seahra \cite{HS} impose that
$V\to 0$ at infinity and then approximate this for numerical
purposes by imposing $V(t,r_c)=0$ on an outer boundary $r=r_c$ at
all $t>0$.\footnote
{In \cite{HS} the outer boundary is defined by $F(t,r_c)=const$ and
thus $r_c$ can vary in time.}
As a result, they have
\begin{equation}
H(t,r_c)=-\frac{1}{(n-1)}\int\limits_{r_0}^{r_c} H^2(t,r)dr \le 0\ ,
\label{ew3.22}
\end{equation}
and $=0$ only when $H(t,r)=0$ for all $r$, so the $r=r_c$ boundary
is necessarily concave in the direction of increasing $r$ at all
$t>0$. This has the potential to generate negative scalar curvature
at $r=r_c$ since
\begin{equation}
R=-2\frac{\partial H}{\partial r} - \frac{n}{(n-1)}H^2
+ (n-1)(n-2)\left ( S/V_{n-1} \right )^{-2/(n-1)}
\ . \label{ew3.23}
\end{equation}
Then Proposition 2.1 and, concomitantly, Proposition 2.4 would no
longer apply.

\subsection{Area radius coordinates}

\noindent This is also not a parabolic system, but can be reduced to
a single parabolic equation. The technicque was used by \cite{OW} to prove
existence and convergence of rotationally symmetric asymptoticall flat flows.
Unfortunately the technique fails for initial data containing a minimal surface.
For this reason, the results of \cite{OW} do not apply when minimal surfaces
are present. In the rotationally symmetric setting, fix $X$ such
that the components of the flow equation in directions tangent to
the orbits of rotational symmetry are trivial; i.e., so that
$\frac{\partial g_{ij}}{\partial t}=0$ for $i,j\neq 1$. This means
that the metric coefficient $B$ above will be time-independent. We
let $f=e^A$ and so we write the metric as
\begin{equation}
ds^2 = f^2(t,r)dr^2+r^2g(S^{n-1},{\rm can})\ , \label{ew3.24}
\end{equation}
and call $r$ the area radius. Then in local coordinates this fixes
\begin{equation}
X=-\frac{R_{22}}{g_{11}\Gamma^1_{22}}\frac{\partial}{\partial r}
=-\frac{R_{33}}{g_{11}\Gamma^1_{33}}\frac{\partial}{\partial r}= \dots \ .
\label{ew3.25}
\end{equation}
Then(\ref{eq1.2}) becomes a constrained system in which the
constraint can be solved, leading to the single parabolic equation
\begin{equation}
\frac{\partial f}{\partial t}=\frac{1}{f^2}\frac{\partial^2
f}{\partial r^2} -\frac{2}{f^3}\left ( \frac{\partial f}{\partial r}
\right )^2 + \left ( \frac{(n-2)}{r}-\frac{1}{rf^2} \right )
\frac{\partial f}{\partial r}-\frac{(n-2)}{r^2f}\left ( f^2 -1 \right )
\ . \label{ew3.26}
\end{equation}
Having a single parabolic equation instead of a system is an
enormous advantage \cite{OW}, and the maximum principle can be used
to show that $f$ remains bounded whenever it is bounded on the
initial data, yielding uniform parabolicity. However, unboundedness
of $f$ on the initial data corresponds precisely to the presence of
a minimal sphere, since the mean curvature of $r=const$ sphere is
given in this system by $H=\frac{(n-1)}{rf}$. (It is shown in
\cite{OW} that no solution of (\ref{ew3.26}) can form a minimal
sphere during the evolution if none is present initially.)

\subsection{Conformal gauge: $n=2$ dimensions}

\noindent Every rotationally symmetric metric is conformally flat.
However, for $n>2$, the conformal class varies throughout the flow.
To see this, set $A=B$ in the metric (\ref{eq3.2}). Then equations
(\ref{eq3.6}, \ref{eq3.7}) reduce to the single equation
\begin{equation}
\frac{\partial A}{\partial t}= e^{-2A} \left \{
\frac{\partial ^2 A}{\partial r^2} + \frac{(n-1)}{r}
\frac{\partial A}{\partial r} +(n-2) \psi'(r) \left [
\frac{\partial A}{\partial r} + \frac{1}{r} \right ] \right \}
\ , \label{eq3.27}
\end{equation}
together with a restriction on $\psi$:
\begin{equation}
0=(n-2) \left [ \left ( \frac{\partial A}{\partial r} \right ) ^2
+\psi''(r) -2\frac{\partial A}{\partial r} \psi'(r) -\frac{1}{r}\psi'(r)
\right ] \ . \label{eq3.28}
\end{equation}
When $n>2$, we obtain from (\ref{eq3.28}) that
\begin{equation}
\psi'(r)=re^{2A(t,r)}\left [ C(t) - \int \frac{e^{-2A(t,r')}}{r'}\left (
\frac{\partial A}{\partial r'} (t,r')\right )^2 dr'\right ] \ , \label{eq3.29}
\end{equation}
for some function $C(t)$. Since the left-hand side is time-independent, so must
be the right. Thus, the system (\ref{eq3.27}, \ref{eq3.28}) is rarely solvable unless $n=2$.
In $n=2$, however, (\ref{eq3.28}) is trivial and then (\ref{eq3.27}) can be solved.

\section{Numerical results}
\setcounter{equation}{0}

\subsection{Initial data}

\noindent We begin with the class of metrics
\begin{equation}
ds^2=\frac{d\rho^2}{1-(\rho_0/\rho)^{\alpha}}
+\rho^2 g(S^{n-1},{\rm can})\ , \label{eq4.1}
\end{equation}
where $g(S^{n-1},{\rm can})$ is the metric with constant unit
sectional curvature on the $(n-1)$-sphere. For $n=3$ and $\alpha\in
[2,1]$, these metrics were used as initial data by \cite{HS}. When
$\alpha=n-2$, this is a static hypersurface in the
Schwarzschild-Tangherlini metric. There is a minimal surface at
$\rho=\rho_0$.

To obtain initial data for the functions $A(t,r)$ and $S(t,r)$ used
in the evolution equations, we first transform to isotropic
coordinates. We obtain
\begin{eqnarray}
ds^2 &=&\beta^2(r) \left ( dr^2+r^2g(S^{n-1},{\rm can}) \right )
\ , \label{eq4.2}\\
\beta(r)&=&\left ( \frac{\rho_0}{2^{1/\alpha}r_0}\right )
\left (1+\frac{r_0^{\alpha}}{2r^{\alpha}} \right )^{2/\alpha}
\ , \label{eq4.3}
\end{eqnarray}
where $r_0$ is a constant of integration arising in the coordinate
transformation. Areas of constant-$r$ spheres are given by
\begin{equation}
V_{n-1}r^{n-1}
\left ( \frac{\rho_0}{2^{1/\alpha}r_0}\right )^{\frac{2}{n-1}}
\left (1+\frac{r_0^{\alpha}}{2r^{\alpha}}
\right )^{\frac{4}{(n-1)\alpha}} \ , \label{eq4.4}
\end{equation}
where $V_{n-1}$ is the volume of an $(n-1)$-sphere of constant unit
sectional curvature. The minimal sphere now lies at
\begin{equation}
r:=r_0/2^{1/\alpha}\ . \label{eq4.5}
\end{equation}
We choose $r_0=2^{1/\alpha}$, so the minimal sphere occurs always at
$r=1$. We then fix $\rho_0=1$, which fixes the
area of the initial minimal sphere to be $V_{n-1}$. Then the initial
data are
\begin{eqnarray}
A(0,r)&=&\log \beta = \frac{2}{\alpha}\log
\left (1+\frac{1}{r^{\alpha}} \right )
-\frac{2}{\alpha} \log 2 \ , \label{eq4.6}\\
S(0,r)&=&V(r)r^{n-1}
\left ( \frac{1}{2}\right )^{\frac{4}{(n-1)\alpha}}
\left (1+\frac{1}{r^{\alpha}}
\right )^{\frac{4}{(n-1)\alpha}}\ .  \label{eq4.7}
\end{eqnarray}

\subsection{The boundary data}

\noindent Our initial data have reflection symmetry in the minimal
sphere, as well as rotational symmetry. The flow equations preserve
isometries, so the subsequent evolution will share these symmetries.

We therefore choose to place a boundary at the location $r=1$ of the
minimal sphere. The idea is then to take the Riemannian double
$(D,g_D)$ (the topological double, with the Riemannian metric $g_D$
induced by pullback of $g(t)$ under the covering map) at any time
$t$ along the evolution. For this, we must keep the mean curvature
zero at $r=1$. Rotational symmetry then ensures the $r=1$ sphere is
totally geodesic, so the double is a smooth manifold. In general, we
do not expect $R(t,1)=0$ at the $r=1$ boundary. Nonetheless,
the arguments of Section 2 apply to $(D,g_D(t))$ provided
that $R(t,r)\to 0$ on the asymptotic ends ($D$, of course, has two
of them).

As we wish to study the noncompact case, we would like to allow $r$
to range through all values $\ge 1$. However, for numerical
purposes, we must either choose a finite cut-off or use a more
sophisticated method (such as conformal compactification, which then
would introduce the complication of dealing with a singular boundary
value problem). For simplicity, we choose a finite cut-off, so $r\in
[1,r_c]$.

We now have two boundaries and a parabolic system of two PDEs, so we
expect four boundary conditions, two at each boundary. Three are
obvious. These are that the mean curvature should vanish at the
$r=1$ boundary and the DeTurck vector field should vanish at both
boundaries (to make it possible to infer the appropriate conclusions
for Ricci flow from our results for the Hamilton-DeTurck flow on
the same bounded manifold). For the remaining condition at the outer
boundary, we would prefer to hold the scalar curvature constant
(preferably zero) there. Instead, we must set a condition that
contains no worse than first derivatives of $A$ and $S$. We choose
to hold the mean curvature equal to a constant $\lambda$ on the
outer boundary. For the case of a single boundary, Cortissoz
\cite{Cortissoz} showed long-time existence of the flow with this
boundary condition, so it would be interesting to see if any vestige
of this result remains for the present case (specifically,
noncollapse of the minimal surface when $\alpha$ is large).

The conditions become, for all $t>0$,
\begin{eqnarray}
\text{At $r=1$:}&&\begin{cases} \frac{\partial S}{\partial r}=0 \\
  \frac{\partial A}{\partial r} + 1 =0
\end{cases}\label{eq4.8}\\
\text{At $r=r_c$:}&&\begin{cases}e^{-A(t,r_c)}\frac{1}{S}
\frac{\partial S}{\partial r}
=\lambda = const\\
\frac{\partial A}{\partial r}=\psi'(r_c)+e^A\lambda
-(n-1)e^{2A} \left ( \frac{V_{n-1}}{S}
\right )^\frac{2}{n-1} r_c^2\left ( \frac{1}{r_c} + \psi'(r_c) \right )
\end{cases}\label{eq4.9}
\end{eqnarray}
We take $\lambda$ to have the value that it would have in the
background metric at $r=r_c$, so
\begin{equation}
\lambda = (n-1)e^{-\psi(r_c)} \left ( \frac{1}{r_c}
+\psi'(r_c) \right ) \ . \label{eq4.10}
\end{equation}

\subsection{The numerical results}

\noindent For the initial data of section 4.1 and $n=3$ dimensions,
the initial minimal surface is always of area $4\pi$, so $\delta=1$
and thus $a^2=10/3$ is the critical value in Proposition 2.4. This
value is achieved for $\alpha=8/3$.

Our first numerical run probes the case of $n=3$ dimensions, with
$\alpha$ taking values below or near the critical value. Figure 1
shows the results for $\alpha$ values $1$, $2$, $2.5$, and $3$. In
every case, the minimal sphere collapses to zero area. Though we do
not show it, we monitor the scalar curvature at the position of the
minimal sphere and observe that it diverges to $+\infty$ at the
collapse time.
\begin{figure}
\centering
\resizebox{10cm}{!}{\includegraphics{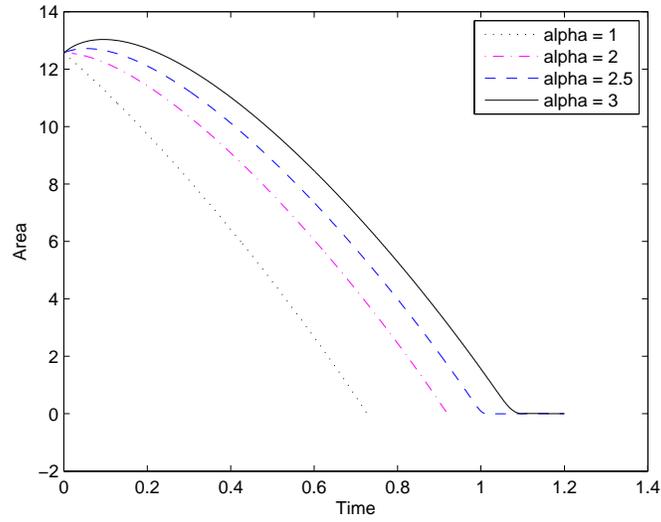}}
\caption{Area of the minimal surface at $r=1$ as a function of
flow time, for $\alpha=1$, $2$, $2.5$, and $3$.}
\end{figure}

In the second run, we probe larger values of $\alpha$, in the hope
of seeing expansion of the minimal surface. The results are
displayed in figure 2. We use $\alpha$ values all the way up to
$\alpha=8$. In each case, collapse eventually commenced, after which
there is no evidence for subsequent re-expansion.
\begin{figure}
\centering
\resizebox{10cm}{!}{\includegraphics{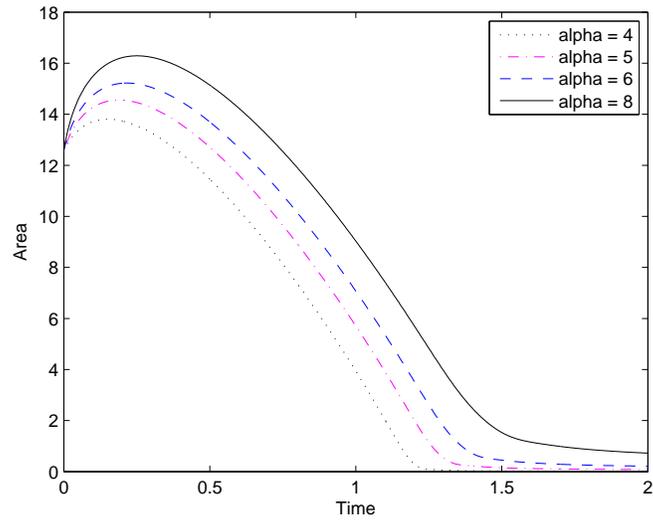}}
\caption{Area of the minimal surface at $r=1$ as a function of
flow time, for $\alpha=4$, $5$, $6$, and $8$.}
\end{figure}
However, at large $\alpha$, the expansion appears eventually to slow
and halt. At such values, as the evolution proceeds, sectional
curvature becomes highly concentrated at the minimal surface,
manifested as a large value of the second spatial derivative of the
area $S$. This may indicate that a singularity forms before collapse
of the minimal surface occurs, but our numerics are not sufficiently
reliable when derivatives become large.

Our third graph deals with the case of $n=2$ dimensions. For initial
data derived from the form of (\ref{eq4.1}) with $\alpha>0$, the
initial scalar curvature will always be negative in 2 dimensions (as
is the case in any dimension when $\alpha>n-2$, so each curve in
figure 3 will initially expand. By (\ref{eq2.20}), we do not have a
collapsing upper barrier function for any $\alpha$, but it is still
possible in principle that collapse could occur. We find no
collapse, despite running for much longer times than for the
higher-dimensional cases.
\begin{figure}
\centering
\resizebox{10cm}{!}{\includegraphics{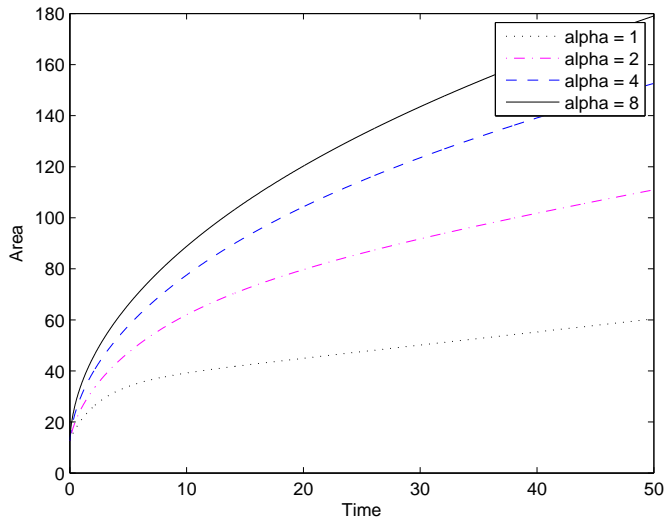}}
\caption{Area of the minimal surface in $n=2$ dimensions as a function of
flow time. The initial data are again given by (\ref{eq4.6}, \ref{eq4.7})
with $n=2$ and various values of $\alpha$.}
\end{figure}

The numerical integrations were performed with MATLAB's ``pdepe''
integration routine for partial differential equations. We performed
several tests of the validity of the code. From equations (\ref{eq3.10}),
(\ref{eq4.6}), and (\ref{eq4.7}), we compute $\frac{\partial S}{\partial t} (0,1)
=4\pi(\alpha-2)$. This serves to define our unit of time, but also serves as a
test of validity of the code, in the sense that the initial derivative should
be linear in alpha and zero at $\alpha=2$. This is verified by figure 4, which
presents a close-up of the evolution near the initial time.
\begin{figure}
\centering
\resizebox{10cm}{!}{\includegraphics{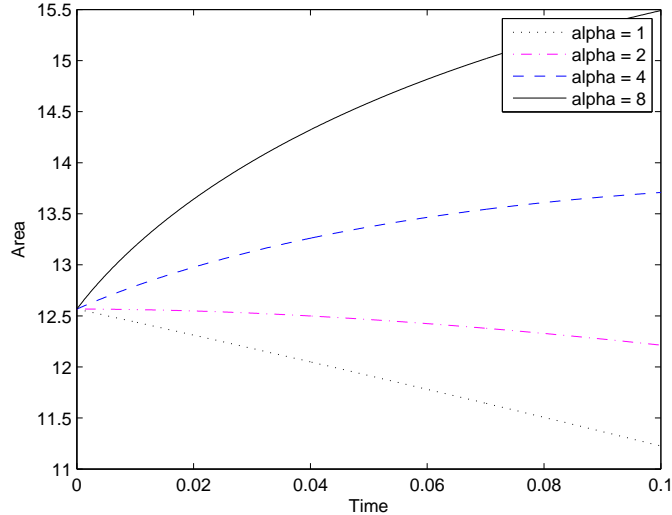}}
\caption{Area of the minimal surface at $r=1$ as a function of
flow time for times close to zero, for $\alpha=1$, $2$, $4$, and $8$.}
\end{figure}
For example, at $\alpha=2$ the slope is zero initially.

As additional verification, we vary the dimension, using initial
data for which the scalar curvature obeys $R(0,r)=0$ in each case.
The initial data may therefore be regarded as the metric on a moment
of time symmetry in an $(n+1)$-dimensional Schwarzschild-Tangherlini
exterior spacetime. We therefore use the evolution equations
(\ref{eq3.6}--\ref{eq3.8}) with arbitrary $n\ge 3$ and set
$\alpha=n-2$ in (\ref{eq4.6}, \ref{eq4.7}). For these runs, we
normalize the initial area of the minimal surface to be
$2\pi^{n/2}/\Gamma(n/2)$. It is easy to verify, using equation
(\ref{eq2.8}), and the fact that $R(0,r)=0$ for such data (and thus
$R(t,r)\ge 0$ by the maximum principle), that for asymptotically
flat boundary conditions the area of the minimal surface will be
bounded above at all times by
\begin{equation}
S(t,0)\equiv\vert \Sigma \vert \le \frac{2\pi^{n/2}}{\Gamma(n/2)} \exp \left [
-\frac12 (n-1)(n-2)t \right ], \ n\ge 3. \label{eq4.11}
\end{equation}
This should therefore give an indication of the effect of our
boundary condition at finite $r=r_c$, which we typically pick to be at
$r_c=10$ (where $r=1$ is the minimal sphere). The results are
depicted in figure 5.
\begin{figure}
\centering
\resizebox{10cm}{!}{\includegraphics{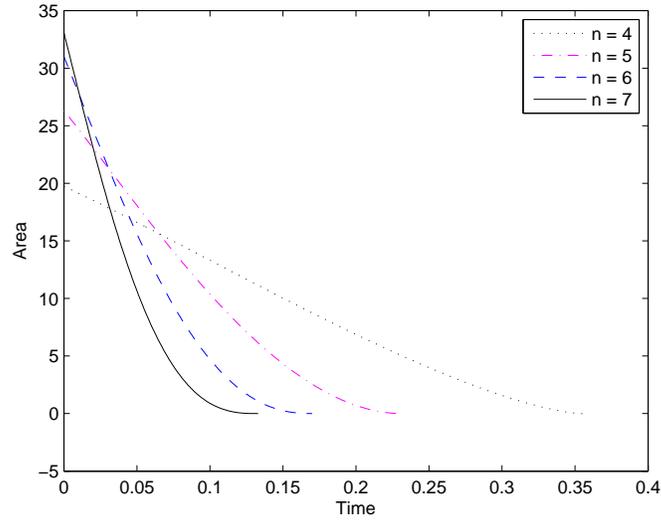}}
\caption{Area of the minimal surface at $r=1$ as a function of
flow time for initial data give by a time-symmetric slice of an
$(n+1)$-dimensional exterior Tangherlini metric for $n=4$, $5$, $6$,
and $7$.}
\end{figure}

A further technical check concerns the choice of different
Hamilton-DeTurck flow. For this, we return to (\ref{eq3.4}) and this
time choose the flat background, so $\psi=0$ in the evolution
equations. The result is displayed in figure 6.
\begin{figure}
\centering
\resizebox{10cm}{!}{\includegraphics{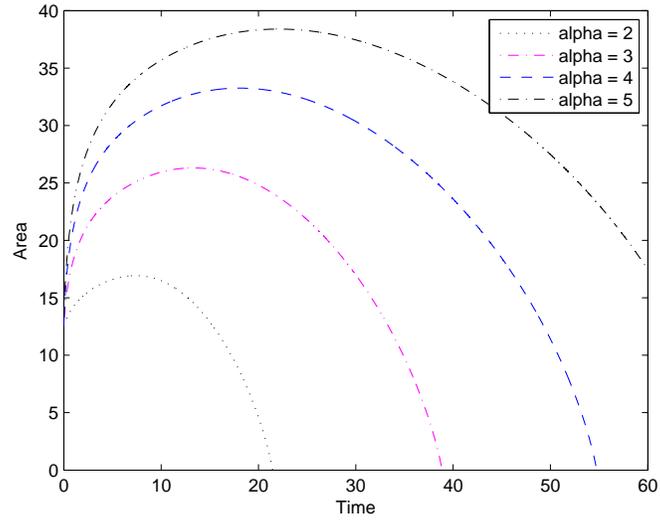}}
\caption{Area of the minimal surface at $r=1$ as a function of
flow time, for $\alpha=2$, $3$, $4$, and $3$. The DeTurck
diffeomorphism is generated by the difference between the flowing connection and
a flat connection.}
\end{figure}
The evolutions take much longer to collapse, but collapse occurs
just as before, and in fact the code seems more reliable close to
the collapse time.\footnote
{The difference appears due to numerical computation of second
spatial derivatives, which grow to much larger magnitude when $\psi$
is the second (nonzero) option in (\ref{eq3.4}) than when $\psi=0$.}
We present graphs for small $\alpha$, but we carried out the
calculation for larger $\alpha$ as well, with similar results.

\section{Discussion}

\noindent The difference between our results and those of
Husain-Seahra may have at least two origins. One is the diffuculty in
dealing numerically with the nonparabolic system (\ref{eq3.16}, \ref{eq3.17}).
Another is the difficulty in
numerically modelling fall-off conditions on a noncompact manifold
using boundary conditions at a fixed boundary. Parabolic equations
exhibit instantaneous action-at-a-distance, though it is typically
exponentially suppressed. Our numerical runs appear to produce some
positive scalar curvature at the large $r$ boundary, which must be
chosen sufficiently distant to suppress the effect on the minimal
surface. The effect is shown in figure 7. For the figure, we
integrated over $r\in [1,100]$ but chose to display only the portion
$r=[90,100]$. At early times we see a small pulse of positive scalar
curvature concentrated at infinity, presumably attributable in some
way to the boundary condition. As time progresses, the pulse
dissipates. The figure shows three times, one shortly after the flow
begins, one at around the midpoint in time of the flow, and one just
before collapse of the minimal surface. We chose to display
$\alpha=3$ but the behaviour is similar for all $\alpha$.
\begin{figure}
\centering
\resizebox{10cm}{!}{\includegraphics{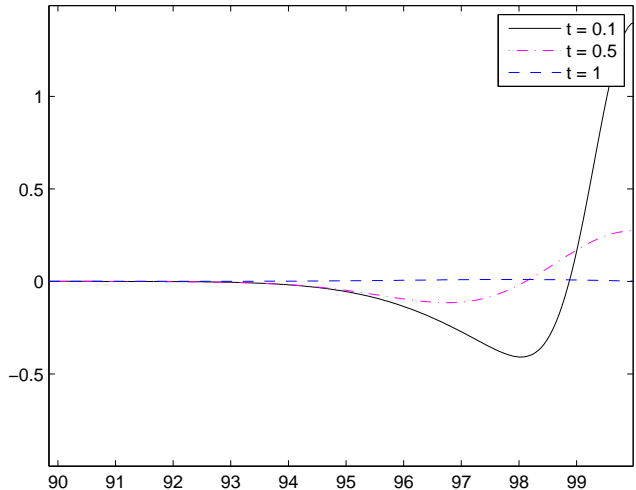}}
\caption{Scalar curvature near the outer boundary at $r=100$ for $\alpha=3$,
plotted at three different times, corresponding roughly to $10\%$, $50\%$,
and $90\%$ of the time to collapse of the minimal surface. At early times,
a disturbance is present near the boundary at $r=100$, but it appears to dissipate
without propagating inward. (For comparison, at $r=1$, the initial scalar curvature
is $R(0,1)=-4$).}
\end{figure}
On the other hand, as discussed in subsection 3.2.2, the boundary
conditions in \cite{HS} imply a concave boundary at $r=r_c$ and,
concomitantly, appear to produce a source of negative curvature
there. Then the maximum principle argument of section 2 would not
apply. The critical behaviour observed in \cite{HS} may be a
property of boundary conditions, rather than evidence for
unstable regions in the space of initial data for asymptotically
flat metrics.

\section{Acknowledgements}

\noindent We are grateful to Paul Mikula, whose numerical code for a
related problem formed the basis from which our code was developed. We thank
Viqar Husain and Sanjeev Seahra for discussions. TB was supported by an
Undergraduate Summer Research Award from the Natural Sciences and Engineering
Research Council of Canada (NSERC). This research was supported by an NSERC
Discovery Grant to EW.

\end{document}